\documentclass[12pt,letterpaper]{article}
\usepackage{amsmath,amsfonts,amsthm,amssymb}

\newcommand\integers{{\mathbb Z}}

\newcommand\reals{{\mathbb R}}

\begin{document}

\newcounter{rom} 

\title{Logic and Partially Ordered Abelian Groups}
\author{David J. Foulis{\footnote{Emeritus Professor, 
University of Massachusetts; foulis@math.umass.edu;
1 Sutton Court, Amherst, MA 01002, USA.}}} 

\date{}
\maketitle

\begin{abstract}

\noindent The unit interval in a partially ordered 
abelian group with order unit forms an interval effect 
algebra (IEA) and can be regarded as an algebraic model 
for the semantics of a formal deductive logic. There 
is a categorical equivalence between the category of 
IEA's and the category of unigroups. In this article, 
we study the IEA-unigroup connection, focusing on 
the cases in which the IEA is a Boolean algebra, an 
MV-algebra, a Heyting MV-algebra, or a quantum logic.  

\end{abstract}

\medskip


\medskip


\section{Introduction}  

Two competing (although not entirely unrelated) 
methods for providing the semantics of a formal deductive 
symbolic calculus ${\cal L}$ are the Kripke many-worlds 
approach and the algebraic-logic approach via 
interpretations of ${\cal L}$ in suitable mathematical 
models. Motivation for the developments in this article 
derives largely from the latter approach via algebraic 
logic.

Typically, a model $L$ for a deductive calculus 
${\cal L}$ is a bounded partially-ordered set equipped 
with operations that qualify it as an algebra (in the 
general sense). An interpretation of ${\cal L}$ 
in $L$ is a mapping $f\mapsto p$ from well-formed 
formulas $f$ of ${\cal L}$ to elements $p$ in $L$ that 
relates the deductive structure of ${\cal L}$ to the 
mathematical structure of $L$. Thus, an element $p\in L$ 
can be considered to be a logical proposition representing 
the equivalence class of all of its antecedent formulas 
$f$. The partial order structure on $L$ is understood in 
the sense that, for propositions $p,q\in L$, $p\leq q$ means 
that $p$ implies $q$ in a manner that is compatible with 
the rules of deduction in ${\cal L}$ for the antecedents of 
$p$ and $q$. Likewise, the algebraic operations on $L$ are  
to be regarded as logical connectives.

Both ${\cal L}$ and its models $L$ are logical systems or 
logical structures, and are often referred to, for short, 
as ``logics."  Thus, in this article, we may refer to a 
partially ordered algebraic structure as a ``logic" when 
we have in mind that it could be construed as a model 
for a deductive symbolic calculus.  

In defining the mathematical structure of a logic $L$, 
one may treat the order and the algebraic structures on 
an equal footing, or one may award primacy to one of 
the structures and derive the other structure from the 
primary one.  For instance, a Boolean algebra $L$ can be 
defined as a bounded, complemented, distributive lattice
or it can be defined as an idempotent ring with unity. 
Starting with the lattice definition, the ring structure 
is derived by taking $x+y$ to be the symmetric difference 
and $xy$ to be the infimum of $x$ and $y$ in $L$.  Starting 
with the ring definition, the order structure is derived 
by defining $x\leq y$ to mean that $x=xy$.

Boolean algebras serve as models for classical propositional 
calculus, poly\-adic Boolean algebras are models for  
first-order predicate calculus, and Heyting algebras are 
models for intuitionistic calculi.  The appropriate models 
for the multi-valued logical calculi of {\L}ukasiewicz 
are the MV-algebras defined in 1957 by C.C. Chang \cite
{CCC}. (See Section 5 below.)
 
In 1986, D. Mundici \cite{Mund} discovered a remarkable 
connection between MV-algebras and lattice-ordered abelian 
groups with order units. (See Section 8 below.) The connection 
is as follows: If $G$ is a lattice ordered abelian group with 
positive cone $G\sp{+}$ and $u\in G\sp{+}$ is an order unit in 
$G$, then the interval $L=G\sp{+}[0,u] =\{e\in G\mid 0\leq 
e\leq u\}$ forms an MV-algebra $(L,0,u,\sp{\perp}\!,\,\hat{+}\,)$, 
where $p\sp{\perp}=u-p$ and the MV-sum is given by $p\,\hat{+}\,
q=(p+q)\wedge u$, for all $p,q\in L$. Conversely, every 
MV-algebra $L$ can be realized as $G\sp{+}[0,u]$ for a 
lattice ordered abelian group $G$ with order unit $u$.  The 
Mundici group $G$ is uniquely determined by $L$ up to 
an isomorphism of partially ordered abelian groups with 
order units.

It turns out that the connection $L\leftrightarrow G$ 
discovered by Mundici admits a considerable generalization 
in which the MV-algebra $L$ is replaced by a so-called 
{\em interval effect algebra} and $G$ is replaced by 
a so-called {\em unigroup}. An interval effect algebra 
may be regarded as a ``logic" in the sense indicated 
above, whence---as per the title of this article---we have 
a connection between a class of logics and a class of 
partially ordered abelian groups.  As some logicians may be 
unfamiliar with the theory of partially ordered abelian 
groups, and some experts on partially ordered abelian 
groups my be uncomfortable with algebraic logic, there is 
a need for an exposition of the logic-unigroup connection.
Our purpose in what follows is to explicate and study 
this connection. Although this article is largely expository, 
a number of unpublished results of N. Ritter \cite{Ritt} 
are cited, and Theorem 8.7, which characterizes the unigroup 
associated with a Heyting effect algebra, is new.  

\section{Boole and the Logic-Algebra Connection} 

In laying the foundations for the algebra that now bears 
his name, George Boole was strongly motivated by analogies 
between the logic of classes and ordinary arithmetic.  For 
classes $x$ and $y$, he took what we now call the 
intersection as the proper interpretation of the ``product" 
$xy$. However, his interpretation of the ``sum" $x+y$ 
differed from what we now call the union in that he 
insisted that it be defined only when $xy=0$. Thus, in his 
1854 masterpiece, {\em The Laws of Thought} \cite{Boo}, Boole 
wrote,

\medskip 

\noindent{\small ``The expression x+y seems indeed 
uninterpretable, unless it be assumed that the things 
represented by x and the things represented by y are 
entirely separate; that they embrace no individuals in 
common."}

\medskip

In {\em The Laws of Thought}, Boole had previously indicated 
that he was well aware of what is now called the union (and even 
the symmetric difference) of classes $x$ and $y$, so his 
decision to write $x+y$ only when $xy=0$ might seem puzzling. 
Indeed, Boole's contemporary, W.S. Jevons, expressed strong 
disagreement with Boole over his unwillingness to give $x+y$ 
an unrestricted interpretation. Using mathematical tools 
unavailable to Boole and Jevons, I. Hailperin has employed   
signed multisets to demonstrate a perfect harmony of Boole's 
product $xy$ and restricted sum $x+y$ with the corresponding 
operations of ordinary arithmetic \cite[pp. 87-112]{IH76}, 
\cite{IH81}. We prefer to recast Hailperin's signed multisets 
in mathematical terms more conducive to the developments in 
this article. To begin with, we assume that the classes 
$x,y,...$ that concerned Boole can be organized into field 
${\cal B}$ of sets. 

\medskip

\noindent{\bf 1.1 Example} Let ${\cal B}$ be a field of 
subsets of a nonempty set $X$ and let $\integers$ be 
the ordered ring of integers. Define ${\cal F}({\cal B},
\integers)$ to be the commutative ring under pointwise 
operations of all bounded functions $f\colon X\to\integers$ 
such that $f\sp{-1}(n)\in\cal{B}$ for all $n\in\integers$.  
The function $1$ that maps all elements in $X$ to the 
integer $1$ is a unity element for the ring ${\cal F}
({\cal B},\integers)$. Under the pointwise partial order, 
${\cal F}({\cal B},\integers)$ is a partially ordered 
(in fact, a lattice-ordered) commutative ring with unity.
The interval $E=\{e\in{\cal F}({\cal B},\integers)
\mid 0\leq e\leq 1\}$ is in bijective correspondence 
with ${\cal B}$ under the mapping $e\leftrightarrow 
e\sp{-1}(1)$. Under this correspondence, the product 
$ef$ in the ring ${\cal F}({\cal B},\integers)$ of 
elements $e,f\in E$ corresponds to the intersection 
$e\sp{-1}(1)\cap f\sp{-1}(1)$ and, if $ef=0$, the 
sum $e+f$ corresponds to the (disjoint) union 
$e\sp{-1}(1)\cup f\sp{-1}(1)$. In this way, the 
restrictions to $E$ of the product and sum in the 
commutative ring ${\cal F}({\cal B},\integers)$ 
match perfectly with Boole's product and sum of 
the corresponding elements of ${\cal B}$.\hspace{\fill}
$\square$

\medskip

By the Stone representation theorem, a Boolean 
algebra $B$ can be represented as the field ${\cal B}$ 
of compact open subsets of a compact, Hausdorff, 
totally-disconnected topological space $X$. Let 
${\cal F}({\cal B},\integers)$ be the partially 
ordered commutative ring with unity in Example 1.1.  
The Boolean algebra $B$ forms an MV-algebra $(B,0,1,
\sp{\perp}\!,\,\hat{+}\,)$, where $x\mapsto x\sp{\perp}$ is 
the Boolean complementation and the MV-sum is given by 
$x\,\hat{+}\,y=x\vee y$. As such, the Mundici group 
corresponding to $B$ is in fact the partially ordered 
additive group of the ring ${\cal F}({\cal B},\integers)$ 
with $1$ as the order unit.

Of course, a Boolean algebra $B$ can be organized into 
a commutative idempotent ring with unity by using 
symmetric difference as the sum and the Boolean meet 
as the product. But, then, $2x=0$ holds for every 
element $x\in B$, a radical departure from ordinary 
arithmetic with which Jevons might have been 
comfortable---but certainly not Boole. On the other hand, 
as in ordinary arithmetic, the additive group of the ring 
${\cal F}({\cal B},\integers)$ is torsion free, i.e., if 
$n$ is a nonzero integer, $f\in{\cal F}({\cal B},\integers)$, 
and $nf=0$, then $f=0$.  

We note that the interval $E$ in Example 1.1 is precisely 
the set of idempotents in the ring ${\cal F}({\cal B},
\integers)$ and thus can be singled out without invoking 
the partial-order structure of the ring.  This would be 
in keeping with the work of Hailperin alluded to above. 
However, it is the partial-order structure that best relates 
to the theme of this article.

\section{Effect Algebras and Boolean Effect Algebras} 

Boole's pioneering work ultimately led to the conception of a 
Boolean algebra, either as a bounded, complemented, distributive 
lattice, or equivalently, as an idempotent ring with unity. 
However, neither of these formulations is based directly on 
Boole's original notions of the product $xy$ and restricted 
sum $x+y$ for classes. 

It is possible to formulate an alternative definition of a 
Boolean algebra involving nothing but the constants $0$, $1$  
and the restricted sum $x+y$. (The product $xy$ then emerges 
as a derived concept.) In formulating this definition, (Definition 
3.5 below), we write $\oplus$ rather than $+$ to emphasize that 
it is only a partially defined operation and also to avoid 
confusing it with addition in the abelian groups to be introduced 
later. Also, we write $u$ instead of $1$ to avoid confusion with 
the numeral $1$. We begin with a basic definition. Apart from the 
change of notation from $+$ to $\oplus$, and from $1$ to $u$, 
axioms (i)--(iv) in Definition 3.1 below are obviously consistent 
with Boole's notion of a restricted sum. 
  
\medskip

\noindent{\bf 3.1 Definition} An {\em effect algebra} is a 
system $(E,0,u,\oplus)$ consisting of a set $E$, special 
elements $0,u\in E$ called the {\em zero} and the {\em unit}, 
and a partially defined binary operation $\oplus$ on $E$, 
satisfying the following conditions for all $x,y,z\in E$:
\begin{list}%
{(\roman{rom})}{\usecounter{rom}
\setlength{\rightmargin}{\leftmargin}} 
\item ({\em Commutativity of $\oplus$}) If $x\oplus y$ is 
 defined, then $y\oplus x$ is defined and $x\oplus y=y\oplus x$. 
\item ({\em Associativity of $\oplus$}) If $x\oplus y$ and 
 $(x\oplus y)\oplus z$ are defined, then $y\oplus z$ and 
 $x\oplus(y\oplus z)$ are defined and $(x\oplus y)\oplus z
 =x\oplus(y\oplus z)$. 
\item ({\em Supplementation property of $\oplus$}) For each 
 $x\in E$, there is a uniquely determined $y\in E$ such that 
 $x\oplus y$ is defined and $x\oplus y=u$. 
\item ({\em Zero-unit property of $\oplus$}) If $x\oplus u$ is 
 defined, then $x=0$. 
\end{list}

Effect algebras were introduced in 1994 \cite {FB} as 
abstractions of the algebra of Hilbert-space effect operators 
used in the study of the theory of measurement in quantum  
mechanics \cite{BLM}.  

\medskip

\noindent{\bf 3.2 Example} Let $R$ be a (not necessarily 
commutative) ring with unity $1$ and let $E$ be the set of 
idempotents in $R$.  If $e,f\in E$, let $e\oplus f :=e+f$ 
iff $ef=fe=0$. (We use := to mean ``equals by definition 
and ``iff" to mean ``if and only if".) Then $(E,0,1,\oplus)$ 
is an effect algebra.\hspace{\fill}$\square$

\medskip

In accord with mathematical tradition, we often say 
that $E$ is an effect algebra when we really mean that 
$(E,0,u,\oplus)$ is an effect algebra. 

\medskip

\noindent{\bf 3.3 Definition} Let $E$ be an effect algebra 
and let $x,y,z\in E$.  Then:
\begin{list}%
{(\roman{rom})}{\usecounter{rom}
\setlength{\rightmargin}{\leftmargin}}
\item We say that $x$ and $y$ are {\em orthogonal} and 
 write $x\perp y$ iff $x\oplus y$ is defined. If $x\perp 
 y$, then $x\oplus y$ is called the {\em orthogonal sum} 
 or for short, the {\em orthosum} of $x$ and $y$.  If we 
 assert that $x\oplus y=z$, we understand that, necessarily, 
 $x\perp y$.
\item If there exists $z\in E$ such that $x\oplus z=y$, we 
 say that $x$ is {\em less than or equal to} $y$ and write 
 $x\leq y$.
\item $x\sp{\perp}$, called the {\em supplement of $x$}, 
 denotes the unique element in $E$ such that $x\oplus 
 x\sp{\perp}=u$.
\end{list}

For Boole's classes $x$ and $y$, the relation $x\perp y$ 
in Definition 3.3 would correspond to the requirement that 
$xy=0$, and the relation $x\leq y$ would correspond to 
the condition that the things represented by $x$ are among 
the things represented by $y$. Also, the ``supplement" 
$x\sp{\perp}$ would be the complement $1-x$ of $x$. 

See \cite{FB} for proofs of the following properties of an 
effect algebra $E$: The relation $\leq$ is a partial order 
relation on $E$ and $x\in E\Rightarrow 0\leq x\leq u$. 
Also, $\leq$ satisfies the following {\em cancellation law}: 
If $x,y,z\in E$, $x\perp z$, and $y\perp z$, then $x\oplus 
z\leq y\oplus z\Rightarrow x\leq y$. Furthermore, 
$x\perp y\Leftrightarrow x\leq y\sp{\perp}$, $x\leq y 
\Rightarrow y\sp{\perp}\leq x\sp{\perp}$, $(x\sp{\perp})
\sp{\perp}=x$, $x\perp 0$, $x\oplus 0=x$, $0\sp{\perp}=
u$, and $u\sp{\perp}=0$. 

If an effect algebra $E$ is regarded as a ``logic" in the sense 
alluded to in Section 1 (i.e., as an algebraic model for a 
deductive logical calculus), then elements $x,y\in E$ can be 
thought of as ``propositions,"  $x\leq y$ means that $x$ 
``implies" $y$, and $0,u\in E$ are ``anti-tautological" and 
``tautological" constants, respectively. The condition $x\perp 
y$ means that, in some sense, the propositions $x$ and $y$ 
``refute" each other. The supplementation mapping $x\mapsto 
x\sp{\perp}$ is a (perhaps attenuated) version of ``logical 
negation," and the ``double negation law" $x=(x\sp{\perp})\sp
{\perp}$ holds. If $x\perp y$, then $x\oplus y$ is to be regarded 
as a sort of (perhaps rarefied) version of ``logical disjunction" 
of the mutually refuting propositions $x$ and $y$. Thus, in 
Definition 3.1, property (iii) may be considered to be a 
rendition of the ``law of the excluded middle" and property (iv) 
may be regarded as a (very) weak ``law of consistency" \cite
[Definition 5.1.1]{DGG}.  Parts (i) and (ii) of the following 
definition are motivated by the observation that, if $x\oplus y$ 
is a logical disjunction of $x$ and $y$ in more or less the 
classical sense, then $x,y\leq p\Rightarrow x\oplus y\leq p$ 
for each proposition $p\in E$.  Part (iii) is suggested by a 
slightly more subtle property of classical logical 
disjunction. 

\medskip

\noindent{\bf 3.4 Definition}  Let $E$ be an effect algebra.    
Then:
\begin{list}%
{(\roman{rom})}{\usecounter{rom}
\setlength{\rightmargin}{\leftmargin}}
\item An element $p\in E$ is {\em principal} iff, for 
 all $x,y\in E$, the conditions $x\perp y$ with $x,y\leq p$ 
 imply that $x\oplus y\leq p$.
\item $E$ is an {\em orthomodular poset} iff every element 
 $p\in E$ is principal.
\item $E$ has the {\em Riesz-decomposition property} 
 iff, for all $x,y,z\in E$, if $y\perp z$ and $x\leq 
 y\oplus z$, there exist $x\sb{1},x\sb{2}\in E$ 
 such that $x\sb{1}\leq y$, $x\sb{2}\leq z$, and 
 $x=x\sb{1}\oplus x\sb{2}$.
 \end{list}

In Definition 3.4 (iii), note that it is not necessary to 
assume that $x\sb{1}\perp x\sb{2}$ since the facts that 
$x\sb{1}\leq y$, $x\sb{2}\leq z$, and $y\perp z$ imply that 
$x\sb{1}\leq y\leq z\sp{\perp}\leq(x\sb{2})\sp{\perp}$, 
whence $x\sb{1}\perp x\sb{2}$. 

Let $(B,\leq,0,u,\wedge,\vee)$ be a Boolean algebra, regarded 
as a bounded distributive lattice, and organize $B$ into an 
effect algebra $(B,0,u,\oplus)$ with $x\oplus y :=x\vee y$ 
iff $x\wedge y=0$, for all $x,y\in B$. Then the effect-algebra 
inequality $\leq$ in Definition 3.3 (ii) coincides with the 
Boolean inequality and the effect-algebra supplement $x\sp
{\perp}$ of $x\in B$ coincides with the Boolean complement of 
$x$. If $x,y,p\in B$, $x\perp y$, and $x,y\leq p$, then $x\oplus 
y=x\vee y\leq p$, so $p$ is principal, and therefore $B$ is an 
orthomodular poset. Also, if $x,y,z\in B$, with $y\perp z$ and 
$x\leq y\oplus z=y\vee z$, then with $x\sb{1} :=x\wedge y\leq y$ 
and $x\sb{2} :=x\wedge z\leq z$, we have $x=x\sb{1}\vee x\sb{2}=
x\sb{1}\oplus x\sb{2}$, so $B$ has the Riesz-decomposition 
property.  Thus, every Boolean algebra is a ``Boolean effect 
algebra" as per the following definition.

\medskip

\noindent{\bf 3.5 Definition} A {\em Boolean effect algebra} 
is an orthomodular poset with the Riesz-decomposition property.

\medskip

\noindent{\bf 3.6 Theorem} {\em As a bounded partially ordered 
set, every Boolean effect algebra is a complemented distributive 
lattice, i.e., a Boolean algebra, in which the supplement of 
each element coincides with its Boolean complement.} 

\medskip

As a consequence of Theorem 3.6 and the remarks preceeding 
Definition 3.5, {\em Boolean algebras are mathematically 
equivalent to Boolean effect algebras}. But notice that, in 
Definitions 3.1, 3.3, 3.4, and 3.5, there is no direct reference 
to the meet and join operations $\wedge$ and $\vee$. The latter 
operations arise from the algebra of $0$, $u$, and $\oplus$, 
rather than vice versa. A proof of Theorem 3.6 will emerge from 
the subsequent developments in this section.

Let $x$ and $y$ be elements of an effect algebra $E$. We write 
the meet (i.e., the infimum, or the greatest lower bound) of 
$x$ and $y$, if it exists in the partially ordered set 
$(E,\leq)$, as $x\wedge y$. Likewise, the join (i.e., the 
supremum, or the least upper bound) of $x$ and $y$, if it 
exists in $(E,\leq)$, is written as $x\vee y$. If we write an 
equation of the form $x\wedge y=z$, we mean that $x\wedge y$ 
exists and equals $z\in E$, and a similar convention holds for 
$x\vee y$. As $x\mapsto x\sp{\perp}$ is order inverting and of 
period two, we have the {\em De Morgan laws}---({\em for meet})  
if $x\wedge y$ exists, then $(x\wedge y)\sp{\perp}=x\sp{\perp}
\vee y\sp{\perp}$, and ({\em for join}) if $x\vee y$ exists, then 
$(x\vee y)\sp{\perp}=x\sp{\perp}\wedge y\sp{\perp}$. (Caution: 
In the general case, if one regards $E$ as a logic, the 
question of whether $x\wedge y$ and $x\vee y$, when they exist, 
should be construed as the conjunction and disjunction of the 
propositions $x$ and $y$ presents subtleties \cite{JP}.)

In the literature, an orthomodular poset is usually defined 
as a structure $(E,\leq,0,u,\sp{\perp})$ consisting of a 
bounded partially ordered set $(E,\leq,0,u)$ together with 
an order-reversing mapping $\sp{\perp}\colon E\to E$ of 
period two such that, (i) for all $x,y\in E$, $x\leq y\sp{\perp}
\Rightarrow x\vee y \text{ exists in }E$, (ii) $x\vee x\sp
{\perp}=u$, and (iii) $x\leq y\Rightarrow y=x\vee(x\vee 
y\sp{\perp})\sp{\perp}$ \cite{Fou62,Kalm,PP}. Condition (iii) 
is called the {\em orthomodular identity}.  To organize such 
a structure into an orthomodular poset according to Definition 
3.4 (ii), one defines $x\oplus y=x\vee y$ iff $x\leq y\sp{\perp}$.
Conversely, it is not difficult to verify that an orthomodular 
poset as per Definition 3.4 is an orthomodular poset according 
to the traditional definition.

\medskip

\noindent{\bf 3.7 Definition} The effect algebra $E$ is 
{\em lattice ordered} iff, as a bounded partially ordered set 
$(E,\leq,0,u)$, it forms a lattice $(E,\leq,0,u,\wedge,\vee)$, 
i.e., $x\wedge y$ and $x\vee y$ exist for all $x,y\in E$. 
If $E$ is a lattice-ordered effect algebra and the lattice 
$(E,\leq,0,u,\wedge,\vee)$ is distributive, we say that 
$E$ is a {\em distributive} effect algebra.

\medskip

Traditionally, an {\em orthomodular lattice} \cite{Beran,Kalm,PP} 
is defined as an orthomodular poset that is also a lattice. 
Thus, from the point of view of effect algebras, {\em 
an orthomodular lattice is a lattice-ordered effect algebra in 
which every element is principal}, and {\em a Boolean effect 
algebra is a distributive orthomodular lattice.}    

A distributive effect algebra is not necessarily a Boolean 
effect algebra.  For instance, the unit interval 
$[0,1]\subseteq\reals$ is organized into a distributive 
effect algebra with $u=1$ by defining $x\oplus y :=x+y$ 
iff $x+y\in [0,1]$ for $x,y\in[0,1]$.  The resulting effect 
algebra has the Riesz-decomposition property, but it is not 
a Boolean effect algebra because the only principal elements 
in $[0,1]$ are $0$ and $1$. In fact, $[0,1]$ is a non-Boolean 
MV-effect algebra (see Section 5 below).

\medskip

\noindent{\bf 3.8 Definition} An element $z$ in an effect algebra 
$E$ is said to be {\em central} in $E$ iff (i) both $z$ and 
$z\sp{\perp}$ are principal in $E$, and (ii) for every $x\in E$ 
there are elements $x\sb{1},x\sb{2}\in E$ such that $x\sb{1}
\leq z$, $x\sb{2}\leq z\sp{\perp}$, and $x=x\sb{1}\oplus x\sb{2}$.
The set of all central elements of $E$ is denoted by $C(E)$ and 
called the {\em center} of $E$ \cite{GFP}.

\medskip

Clearly, if the effect algebra $E$ has the Riesz-decomposition 
property, then condition (ii) in Definition 3.8 holds automatically.
Also, every element $z$ in an orthomodular poset satisfies condition 
(i) in Definition 3.8. Consequently, a Boolean effect algebra $E$ is 
its own center, i.e., $C(E)=E$.

\medskip

\noindent{\bf 3.9 Definition} If $E$ is an effect algebra, then 
a subset $S\subseteq E$ is called a {\em subeffect algebra} of 
$E$ iff $0,u\in S$, $x\in S\Rightarrow x\sp{\perp}\in S$, and 
for all $x,y\in S$, $x\perp y\Rightarrow x\oplus y\in S$.

\medskip

If $S$ is a subeffect algebra of the effect algebra $E$, then 
$S$ forms an effect algebra in its own right under the 
restriction to $S$ of $\oplus$.  By \cite{GFP}, the center 
$C(E)$ of an effect algebra $E$ is a subeffect algebra of $E$ 
and $C(E)$ is a Boolean algebra (hence a Boolean effect algebra).
The promised proof of Theorem 3.6 is now at hand, since if $E$ 
is a Boolean effect algebra, then $C(E)=E$, whence $E$ is a 
Boolean algebra.

\medskip

An alternative characterization of Boolean effect algebras can 
be formulated in terms of the notion of compatibility in the 
next definition.  

\medskip

\noindent{\bf 3.10 Definition} Let $E$ be an effect algebra. We say 
that $x,y,z\in E$ are {\em jointly orthogonal} iff $x\perp y$ and 
$(x\oplus y)\perp z$, (whence $y\perp z$ and $x\perp(y\oplus z)$).  
If $x,y\in E$, then $x$ and $y$ are said to be {\em compatible} (or, 
{\em Mackey compatible}), in symbols $xCy$, iff there are jointly 
orthogonal elements $x\sb{1},y\sb{1},z\in E$ such that $x=x\sb{1}
\oplus z$ and $y=y\sb{1}\oplus z$.

\medskip

If $E$ is an orthomodular poset, then $C(E)=\{z\in E\mid zCx\text
{ for all }x\in E\}$, hence {\em a Boolean effect algebra is the 
same thing as an orthomodular poset in which every pair of elements 
is compatible.} 

\medskip

The complete title of Boole's 1854 classic is {\em An 
Investigation of the Laws of Thought on which are founded 
the Mathematical Theories of Logic and Probabilities}. We note 
that Boole's restricted sum interacts perfectly with probability 
assignments $p$ in that $p(x+y)=p(x)+p(y)$ holds for $xy=0$. 
This leads us to the following definition.

\medskip

\noindent{\bf 3.11 Definition}  Let $E$ be an effect algebra. 
If $K$ is an additive abelian group, then a {\em $K$-valued 
measure} on $E$ is a mapping $\phi\colon E\to K$ such that, 
for all $x,y\in E$, $x\perp y\Rightarrow\phi(x\oplus y)=
\phi(x)+\phi(y)$.  Regarding the ordered field $\reals$ of 
real numbers as an additive abelian group, we define a 
{\em probability measure} on $E$ to be an $\reals$-valued 
measure $\pi\colon E\to\reals$ that is {\em positive} 
in the sense that $0\leq\pi(x)$ for all $x\in E$ and 
{\em normalized} in the sense that $\pi(u)=1$.  Denote by 
$\Pi(E)$ the set of all probability measures on $E$. A 
subset $\Delta\subseteq\Pi(E)$ is {\em order-determining} 
iff, for $x,y\in E$, the condition $\pi(x)\leq \pi(y)$ for 
every $\pi\in\Delta$ implies that $x\leq y$.

\medskip

The set $\Pi(E)$ is a convex subset of the real vector space 
under pointwise operations of all mappings $\rho\colon E\to 
\reals$. We denote by $\partial\sb{e}\Pi(E)$ the set of all 
extreme points of $\Pi(E)$. If the elements of $E$ are regarded 
as ``propositions," then a probability measure $\pi\in\Pi(E)$, 
and especially a $\pi\in\partial\sb{e}\Pi(E)$, can be regarded 
as a (possibly multi-valued) ``truth combination" assigning a 
``truth value" $\pi(x)$ on a scale from $0$ (false) to $1$ (true) 
for each proposition $x\in E$. 

If $B$ is a Boolean effect algebra, then $\partial\sb{e}\Pi(B)$ 
is order-determining, elements of $\partial\sb{e}\Pi(B)$ are 
$\{0,1\}$-valued, and $\partial\sb{e}\Pi(B)$ may be identified 
with the Stone space of $B$. That $\partial\sb{e}\Pi(B)$ is order 
determining accounts for the fact that truth tables provide an    
algorithmic decision procedure for classical propositional calculus.   

\section{Quantum Logics} 

Certain effect algebras $E$ can be considered to be algebraic 
models for the semantics of the ``quantum logics," that arise in 
the study of reasoning in quantum theory \cite{DGG}. Rather than 
saying that such an $E$ is an algebraic model for a quantum logic, 
we shall say, for short, that $E$ is a quantum logic.

The genesis of quantum logic was von Neumann's observation 
\cite[p. 253]{JvN},

\medskip

\noindent{\small ``.. the relation between the properties 
of a physical system on the one hand, and the projections 
on the other, makes possible a sort of logical calculus with 
these."}

\medskip

The projections to which von Neumann referred are the bounded 
self-adjoint idempotent operators $P=P\sp{\ast}=P\sp{2}$ on a 
Hilbert space ${\cal H}$, and these projections band together to 
form an orthomodular lattice $\mathbb{P}({\cal H})$. In quantum 
mechanics, the question of whether two projections $P$ and $Q$ 
commute, i.e., whether $PQ=QP$, is of considerable significance, 
and it is important to note that it can be settled strictly in 
terms the structure of $\mathbb{P}({\cal H})$ as an effect 
algebra. In fact, two projections $P$ and $Q$ on the Hilbert 
space ${\cal H}$ are compatible in the orthomodular lattice 
$\mathbb{P}({\cal H})$ iff $PQ=QP$.

Von Neumann's observation that $\mathbb{P}({\cal H})$ can be 
regarded as a logical calculus led to the study of more general 
orthomodular lattices as possible quantum logics. Indeed, S. 
Gudder and others were able to show that much of the theory of 
spectral measures and quantum probability carries over to the 
more general context of an orthomodular lattice that admits 
sufficiently many probability measures \cite{Gud65}. However, 
difficulties associated with the interpretation (as logical 
connectives) of the meet and join of noncommuting projections 
subsequently led to the consideration of more general 
orthomodular posets as quantum logics \cite{Gud68}. Further 
difficulties arising from the necessity of dealing with 
coupled quantum-mechanical systems led to the study of 
orthoalgebras as quantum logics \cite{Fou89}. 

An {\em orthoalgebra} is an effect algebra $E$ such that $x\wedge 
x\sp{\perp}=0$ for all $x\in E$. Every orthomodular poset is an 
orthoalgebra, but not vice versa. By the De Morgan law, every element 
$x$ in an orthoalgebra $E$ satisfies both $x\wedge x\sp{\perp}=0$ and 
$x\vee x\sp{\perp}=u$, i.e., just as in a Boolean algebra, $x
\sp{\perp}$ is a complement of $x$ in $E$.  Thus, in an orthoalgebra, 
we have a semantic version $x\wedge x\sp{\perp}=0$ of the classical 
law of noncontradiction (ex contradictione quodlibet, or Duns Scotus' 
law) and also the excluded middle law (tertium non datur) $x\vee x\sp
{\perp}=u$. If $E$ is an orthoalgebra and $x,y\in E$, then $xCy$ iff 
there is a Boolean subeffect algebra $B$ of $E$ such that $x,y\in B$.  
Also, for an orthoalgebra $E$, the center $C(E)$ is given by $C(E)=
\{z\in E\mid zCx\text{ for all }x\in E\}$.   

In the contemporary theory of quantum measurement \cite{BLM} the 
projection-valued measures favored by von Neumann are replaced by 
more general measures defined on a $\sigma$-field of sets and 
taking on values in the set $\mathbb{E}({\cal H})$ of effect 
operators on a Hilbert space ${\cal H}$. An {\em effect operator} 
on ${\cal H}$ is a bounded self-adjoint operator $A$ on $H$ such 
that $\mathbf{0}\leq A\leq\mathbf{1}$, and the set $\mathbb{E}
({\cal H})$ can be organized into an effect algebra 
$(\mathbb{E}({\cal H}),\mathbf{0},\mathbf{1},\oplus)$, where, for 
$A,B\in\mathbb{E}({\cal H})$, $A\oplus B :=A+B$ iff $A+B\leq
\mathbf{1}$. As such, the effect-algebra partial order coincides 
with the restriction to $\mathbb{E}({\cal H})$ of the usual partial 
order on bounded self-adjoint operators, and if $A\in\mathbb{E}
({\cal H})$, then $A\sp{\perp}=\mathbf{1}-A$. If $\langle\cdot,
\cdot\rangle$ is the inner product on ${\cal H}$, then each unit 
vector $\psi\in{\cal H}$ determines a probability measure $\pi\sb
{\psi}\in\partial\sb{e}\Pi(\mathbb{E}({\cal H}))$ according to 
$\pi\sb{\psi}(A) :=\langle A\psi,\psi\rangle$ for every $A\in
\mathbb{E}({\cal H})$. Therefore, the effect algebra $\mathbb{E}
({\cal H})$ carries an order-determining set of probability measures.  

The orthomodular lattice $\mathbb{P}({\cal H})$ of projection 
operators on ${\cal H}$ is a subeffect algebra of $\mathbb{E}
({\cal H})$, and if $P\in\mathbb{E}({\cal H})$, then $P\in\mathbb{P}
({\cal H})\Leftrightarrow P\text{ is principal}\Leftrightarrow 
P\wedge P\sp{\perp}=\mathbf{0}$. In the passage from the 
orthomodular lattice $\mathbb{P}({\cal H})$ to the larger effect 
algebra $\mathbb{E}({\cal H})$, the supplementation mapping 
$A\mapsto A\sp{\perp}=\mathbf{1}-A$ loses its character as a 
complementation, and $\mathbb{E}({\cal H})$ becomes an algebraic 
model for a {\em paraconsistent logic} \cite{DG89}.

The Hilbert-space effect algebra $\mathbb{E}({\cal H})$ is the 
prototypic effect algebra; both it and its subeffect algebra 
$\mathbb{P}({\cal H})$ are the prototypic quantum logics. Nowadays, 
an effect algebra $E$ is regarded as a quantum logic only if it 
satisfies some of the special properties of the prototypes 
$\mathbb{E}({\cal H})$ or $\mathbb{P}({\cal H})$. Paramount among 
these properties are conditions relating to probability measures, 
especially the condition that $\Pi(E)$ is order determining.  
We shall resist the temptation to give a formal definition 
of a quantum logic. (For authoritative literature on the question of 
just what constitutes a quantum logic, see \cite{DGG} and \cite{PP}.)
However, if ${\cal A}$ is a unital C$\sp{\ast}$-algebra, we propose 
to regard $E :=\{e\in{\cal A}\mid e=e\sp{\ast}\text{ and }0\leq 
e\leq 1\}$, as well as its subeffect algebra $P :=\{p\in{\cal A}
\mid p=p\sp{\ast}=p\sp{2}\}$ as bona fide quantum logics.

\section{MV-Algebras} 

Material in this section is adopted from \cite{FouMV}. The 
following definition is based on \cite[Lemma 2.6]{Mund}.

\medskip

\noindent{\bf 5.1 Definition}  An {\em MV-algebra} is a system 
$(E,0,u,\sp{\perp}\!,\,\hat{+}\,)$ consisting of a set $E$, special 
elements $0,u\in E$ called the {\em zero} and the {\em unit}, a 
unary operation $p\mapsto p\sp{\perp}$ called {\em supplementation} 
on $E$, and a binary operation $\;\,\hat{+}\,\;$ called the {\em 
MV-sum} on $E$ that satisfies the following axioms for all $p,q,r
\in E$:
\[
\begin{array}{lll}
\text{(i) } p\,\hat{+}\,(q\,\hat{+}\,r)=(p\,\hat{+}\,q)\,\hat{+}\,r & 
\text{(ii) } p\,\hat{+}\,q=q\,\hat{+}\,p & \text{(iii) }p\,\hat{+}
\,0=p\\
\text{(iv) }p\,\hat{+}\,u=u & \text{(v) }p\sp{\perp\perp}=p & 
\text{(vi) }0\sp{\perp}=u\\
\text{(vii) }p\,\hat{+}\,p\sp{\perp}=u\\
\text{(viii) }(p\,\hat{+}\,q
 \sp{\perp})\sp{\perp}\,\hat{+}\,p=(q\,\hat{+}\,p\sp{\perp})\sp{\perp}
 \,\hat{+}\,q.
\end{array}
\]

\medskip

\noindent{\bf 5.2 Definition} An {\em MV-effect algebra} is a 
lattice-ordered effect algebra with the Riesz-decomposition 
property.

\medskip

According to the following theorem, originally proved by Chovanec 
and K\^{o}pka \cite{CK} and here translated into the language of 
effect algebras, MV-algebras and MV-effect algebras are 
mathematically equivalent notions.

\medskip

\noindent{\bf 5.3 Theorem} {\em An MV-algebra $(E,0,u,\sp{\perp}\!,
\,\hat{+}\,)$ forms an MV-effect algebra $(E,0,u,\oplus)$ where, for 
$p,q\in E$, $p\oplus q :=p\,\hat{+}\,q$ iff $p\leq q\sp{\perp}$. 
Moreover, for $p,q\in E$, $p\vee q=(p\,\hat{+}\,q \sp{\perp})\sp
{\perp}\,\hat{+}\,p$.  Conversely, an MV-effect algebra $(E,0,u,
\oplus)$ forms an MV-algebra $(E,0,u,\sp{\perp}\!,\,\hat{+}\,)$ 
where $p\mapsto p\sp{\perp}$ is the effect-algebra 
supplementation map and, for $p,q\in E$, $p\,\hat{+}\,q :=
p\oplus(p\sp{\perp}\wedge q)$.}

\medskip

\noindent{\bf 5.4 Corollary} {\em An MV-effect algebra is 
a distributive effect algebra.}

\medskip

\noindent{\bf 5.5 Theorem} {\em If $E$ is a lattice-ordered effect 
algebra, then $E$ is an MV-effect algebra iff, for all $p,q\in E$, 
$p\wedge q=0\Rightarrow p\perp q$.}

\medskip

\noindent{\em Proof}  See \cite[Theorem 3.11]{BFPhi}.\hspace{\fill}
$\square$ 

\medskip

Clearly, every Boolean effect algebra is an MV-effect algebra. 
Furthermore, every MV-effect algebra is an extension of a 
Boolean subeffect algebra, namely its center.

\medskip

\noindent{\bf 5.6 Theorem} {\em Let $E$ be an MV-effect algebra. 
Then}\[
 C(E)=\{c\in E\mid c\wedge c\sp{\perp}=0\}=\{c\in E\mid 
 c\,\hat{+}\,c=c\}.
\]

\medskip

\noindent{\em Proof} See \cite[Theorem 6.1]{FouMV}.\hspace{\fill}
$\square$ 

\section{Heyting and Heyting Effect Algebras} 
 
\noindent{\bf 6.1 Definition} A {\em Heyting algebra} is a system 
$(H,\leq,0,1,\wedge,\vee,\supset)$ such that $(H,\leq,0,1,\wedge,
\vee)$ is a bounded lattice and $\supset$ is a binary operation 
on $H$, called the {\em Heyting conditional}, such that for all 
$p,q,r\in H$, $p\wedge q\leq r\Leftrightarrow p\leq(q\supset r)$.
If $H$ is a Heyting algebra and $p\in H$, then $p\,' :=
(p\supset0)$ is called the {\em Heyting negation} of $p$. A 
Heyting algebra $H$ is called a {\em Stone-Heyting algebra} iff, 
for all $p\in H$, $p\,'\vee (p\,')\,'=1$. 

\medskip

Every Boolean algebra is a Stone-Heyting algebra with the 
material conditional $p\supset q :=p\sp{\perp}\vee q$ as the 
Heyting conditional and $p\,'=p\sp{\perp}$ as the Heyting 
negation. If a Heyting algebra $H$ is a model for the 
semantics of an intuitionistic logic, and if $p,q\in H$, then 
$p\supset q$ is supposed to be a proposition in $H$ asserting 
that $p$ implies $q$. In this regard, we note that $(p
\supset q)=1\Leftrightarrow p\leq q$.  

Let $H$ be a Heyting algebra. Then $(H,\leq,0,1,\wedge,
\vee)$ is a bounded distributive lattice. Also, the Heyting 
negation mapping $'\colon H\to H$ satisfies $p\wedge 
q=0\Leftrightarrow q\leq p\,'$ for all $p,q\in H$. In 
particular, $p\wedge p\,'=0$, so the Heyting negation 
satisfies Duns Scotus' law.  However, tertium non datur 
does not necessarily hold, i.e., $p\vee p\,'=1$ may fail. 
Also, although $p\leq p\,'' :=(p\,')\,'$ always holds, the 
condition $p\,''\leq p$ may fail. In fact, the set $\{c\in 
H\mid c=c\,''\}$, which is the same as the set $H\,' :=\{p\,'
\mid p\in H\}$, forms a Boolean algebra under the restriction 
of the partial order on $H$.  If $p,q\in H\,'$, then $p\wedge 
q\in H\,'$ is the infimum of $p$ and $q$ in $H\,'$; however, 
unless $H$ is a Stone-Heyting algebra, $p\vee q$ need not 
belong to $H\,'$, and the supremum of $p$ and $q$ in $H\,'$ 
is $(p\vee q)\,''=(p\,'\wedge q\,')\,'$. The restriction to 
$H\,'$ of the Heyting negation is the Boolean complementation 
on $H\,'$. We call the Boolean algebra $H\,'$ the {\em Heyting 
center} of $H$.

\medskip

\noindent{\bf 6.2 Definition} A {\em Heyting effect algebra} 
is a lattice ordered effect algebra $E$ equipped with a 
binary operation $\supset$ such that $(E,\leq,0,u,\wedge,
\vee,\supset)$ is a Heyting algebra.

\medskip

Although Boolean algebras are coextensive with Boolean effect 
algebras and MV-algebras are coextensive with MV-effect algebras, 
there are Heyting algebras, and even Stone-Heyting algebras, that 
cannot be organized into Heyting effect algebras.

\medskip

\noindent{\bf 6.3 Theorem} {\em Let $E$ be a Heyting effect 
algebra and let $e,f\in E$. Then:} (i) $e\,'=(e\supset 0)
\in C(E)$. (ii) $e\,'\leq e\sp{\perp}$ {\em with equality 
iff $e\in C(E)$.} (iii) {\em The Heyting center of $E$ 
coincides with the effect-algebra center $C(E)$.} 
(iv) $e\wedge f=0\Rightarrow e\perp f$. (v) {\em $E$ is an 
MV-effect algebra.} (vi) {\em $E$ is a Stone-Heyting algebra.}

\medskip

\noindent{\em Proof} Part (i) follows from \cite[Theorem 3.31]
{Ritt}, and in view of (i), parts (ii)--(vi) follow from 
\cite[Theorems 8.3 and 8.5]{FouMV}.\hspace{\fill}$\square$ 

\medskip 

Since every Heyting effect algebra is an MV-effect algebra, 
we often refer to a Heyting effect algebra as an {\em 
HMV-effect algebra}, or simply as an {\em HMV-algebra}. 
Theorem 6.3 (vi) implies that, unlike Heyting algebras in 
general, the Heyting center of an HMV-algebra is closed under 
the formation of suprema. By \cite[Corollary 3.34]{Ritt}, {\em 
every MV-algebra that is complete as a lattice is an 
HMV-algebra.}  By \cite[Theorem 8.3]{FouMV}, {\em a Heyting 
effect algebra is the same thing as a lattice-ordered effect 
algebra $E$ equipped with a mapping $'\colon E\to C(E)$ such 
that, for all $e,f\in E$, $e\wedge f=0\Leftrightarrow e\leq 
f\,'$.}

\section{Interval Effect Algebras and Unigroups} 

An (additively-written) abelian group $G$ is called a {\em 
partially ordered abelian group} iff it is equipped with a 
partial order $\leq$ that is {\em translation invariant} in 
the sense that, for $g,h,k\in G$, $g\leq h\Rightarrow g+k\leq 
h+k$.  If $G$ is a partially ordered abelian group, then the 
subset $G\sp{+} :=\{g\in G\mid 0\leq g\}$ is called the 
{\em positive cone} in $G$ (in spite of the fact that $0\in
G\sp{+}$).  The positive cone $G\sp{+}$ satisfies the conditions 
(i) $0\in G\sp{+}$, (ii) $g,h\in G\sp{+}\Rightarrow g+h\in G
\sp{+}$, and (iii) $g,-g\in G\sp{+}\Rightarrow g=0$.  Conversely, 
if $G$ is an abelian group and $G\sp{+}\subseteq G$ is a subset 
of $G$ satisfying conditions (i), (ii), and (iii), there is one 
and only one translation-invariant partial order $\leq$ on $G$ 
for which $G\sp{+}$ is the corresponding positive cone, and it 
is determined by $g\leq h\Leftrightarrow h-g\in G\sp{+}$ for all 
$g,h\in G$. A partially-ordered abelian group $G$ is said to be 
{\em archimedean} iff, for all $g,h\in G$, the condition 
$ng\leq h$ for all positive integers $n$ implies that 
$-g\in G\sp{+}$.  

If $G$ is a partially ordered abelian group and $u\in G\sp{+}$, 
define the {\em $u$-interval} $G\sp{+}[0,u] :=\{e\in G\mid 
0\leq e\leq u\}$.  Such a $u$-interval can be organized into 
an effect algebra $(G\sp{+}[0,u],0,u,\oplus)$ by defining 
$p\oplus q :=p+q$ iff $p+q\leq u$, for all $p,q\in G\sp{+}[0,u]$.
As such, the effect-algebra partial order is the restriction to 
$G\sp{+}[0,u]$ of the partial order on $G$, and for $p\in G
\sp{+}[0,u]$, $p\sp{\perp}=u-p$.

A {\em morphism} from an effect algebra $E$ with unit 
$u$ into an effect algebra $F$ with unit $v$ is a mapping 
$\phi\colon E\to F$ such that: (i) $p,q\in E$ with $p\perp q$ 
implies that $\phi(p)\perp\phi(q)$ and $\phi(p\oplus q)=
\phi(p)\oplus\phi(q)$, and (ii) $\phi(u)=v$.  An {\em 
isomorphism} is a bijective morphism $\phi\colon E\to F$ 
such that $\phi\sp{-1}\colon F\to E$ is also a morphism.
If there is an isomorphism $\phi\colon E\to F$, we say that 
$E$ and $F$ are {\em isomorphic}.

\medskip

\noindent{\bf 7.1 Definition} An effect algebra $E$ is called 
an {\em interval effect algebra} (IEA) iff it can be realized as, 
or is isomorphic to, a $u$-interval $G\sp{+}[0,u]$ in a partially 
ordered abelian group $G$ with $u\in G\sp{+}$.

\medskip

By Example 1.1 and the Stone representation theorem, every 
Boolean algebra is an IEA. By Mundici's theorem, every MV-algebra 
is an IEA. Because a Heyting effect algebra is an MV-algebra, 
every Heyting effect algebra is an IEA.  K. Ravindran \cite{Rav96} 
has generalized Mundici's theorem by proving that every effect 
algebra with the Riesz-decomposition property is an IEA. If 
${\cal H}$ is a Hilbert space and $\mathbb{G}({\cal H})$, partially 
ordered in the usual way, is the additive group of bounded 
self-adjoint operators on ${\cal H}$, then by definition 
$\mathbb{E}({\cal H})=\mathbb{G}({\cal H})\sp{+}[\mathbf{0},
\mathbf{1}]$, so the quantum logic $\mathbb{E}({\cal H})$ 
is an interval effect algebra.

Currently, it is not known how to give an intrinsic 
characterization of an IEA.  However, by \cite[Corollary 2.5]
{BFInt}, a subeffect algebra of an IEA is again an IEA, and 
by \cite[Theorem 5.4]{BFInt}, an effect algebra with an 
order-determining set of probability measures is in IEA. As 
a partial converse, it turns out that every IEA admits at 
least one probability measure \cite[Theorem 5.5]{BFInt}.

Let $G$ be a partially-ordered abelian group.  If $G\sp{+}$ 
generates $G$ as a group, then $G$ is said to be {\em directed}.  
It is easy to see that $G$ is directed iff $G=G\sp{+}-G\sp{+}$, 
i.e., iff every element $g\in G$ can be written as $g=g\sb{1}-
g\sb{2}$ with $g\sb{1},g\sb{2}\in G\sp{+}$. If, as a partially 
ordered set, $G$ forms a lattice, then $G$ is called a {\em 
lattice-ordered} abelian group. We say that $G$ has the {\em 
interpolation property} iff, given $a,b,c,d\in G$ with 
$a\leq c$, $a\leq d$, $b\leq c$, and $b\leq d$, there exists 
$t\in G$ such that $a\leq t$, $b\leq t$, $t\leq c$, and 
$t\leq d$ \cite{Good}. A partially-ordered abelian group with 
the interpolation property is called an {\em interpolation 
group}. If $G$ is lattice ordered, it is an interpolation 
group \cite[p. 23]{Good}. 

Suppose that $G$ is a lattice-ordered abelian group and that
$g,h,k\in G$. Then $(G,\wedge,\vee)$ is a distributive lattice, 
$-(g\wedge h)=(-g)\vee(-h)$, $-(g\vee h)=(-g)\wedge(-h)$, 
$(g\wedge h)+k=(g+k)\wedge(h+k)$, $(g\vee h)+k=(g+k)\vee(h+k)$, 
and $g+h=(g\vee h)+(g\wedge h)$ \cite[Chapter 1]{Good}. Define 
$g\sp{+} :=g\vee 0=g-(g\wedge 0)$ and $g\sp{-} :=(-g)\sp{+}=
(-g)\vee 0=-(g\wedge 0)$. Then $0\leq g\sp{+},g\sp{-}$ and $g=
g\sp{+}-g\sp{-}$, hence $G$ is directed. Furthermore, $g\sp{+}
\wedge g\sp{-}=(g+g\sp{-})\wedge g\sp{-}=(g\wedge 0)+g\sp{-}=
(g\wedge 0)-(g\wedge 0)=0$. 

Let $u\in G\sp{+}$.  We say that $u$ is an {\em order unit} 
for $G$ iff, for every $g\in G$, there is a positive integer 
$n$ such that $g\leq nu$. If every $g\in G\sp{+}$ can be 
written as a finite linear combination with positive integer 
coefficients of elements in the $u$-interval $G\sp{+}[0,u]$, 
i.e., if $G\sp{+}[0,u]$ generates $G\sp{+}$ as a semigroup, 
then $u$ is said to be {\em generative} \cite[Definition 3.2]
{BFInt}. If $G$ admits an order unit, then $G$ is directed 
\cite[p. 4]{Good}. If $u$ is generative and $G$ is directed, 
then $u$ is an order unit for $G$ \cite[Lemma 3.1]{BFInt}. 
As a consequence of \cite[Proposition 2.2 (b)]{Good}, if $G$ 
is an interpolation group, and $u$ is an order unit for $G$,  
then $u$ is generative.

\medskip

\noindent{\bf Definition 7.2} A {\em unital group} is a 
partially-ordered abelian group $G$ with a distinguished 
generative order unit $u\in G\sp{+}$, called the {\em unit}.  
If $G$ is a unital group with unit $u$, then the $u$-interval 
$E :=G\sp{+}[0,u]$, regarded as an IEA, is called the {\em 
unit interval} in $G$. A {\em unital homomorphism} from a 
unital group $G$ with unit $u$ into a unital group $H$ with 
unit $v$ is a group homomorphism $\phi\colon G\to H$ such that 
$\phi(G\sp{+})\subseteq H\sp{+}$ and $\phi(u)=v$. A {\em 
unital isomorphism} from $G$ onto $H$ is a bijective unital 
homomorphism $\phi\colon G\to H$ such that $\phi\sp{-1}\colon 
H\to G$ is also a unital homomorphism. Two unital groups $G$ 
and $H$ are {\em isomorphic as unital groups} iff there is a 
unital isomorphism $\phi\colon G\to H$. 

\medskip

Let $E$ be the unit interval in a unital group $G$. Then $E$ 
generates $G\sp{+}$ as a semigroup, and (as $G$ is necessarily 
directed) $G\sp{+}$ generates $G$ as a group; hence, $E$ 
generates $G$ as a group.  Therefore, if $K$ is an abelian 
group, $\Phi\colon G\to K$ is a group homomorphism, and 
$\phi :=\Phi|\sb{E}$ is the restriction of $\Phi$ to $E$, 
then $\Phi$ is uniquely determined by the $K$-valued measure 
$\phi\colon E\to K$. 
      
\medskip

\noindent{\bf 7.3 Definition} Let $G$ be a unital group with 
unit interval $E$. If $K$ is an abelian group, we say that 
$G$ is {\em $K$-universal} iff every $K$-valued measure 
$\phi\colon E\to K$ can be extended to a (necessarily unique) 
group homomorphism $\Phi\colon G\to K$. If $G$ is $K$-universal 
for every abelian group $K$, then  $G$ is called a {\em unigroup} 
\cite{FGB98}.

\medskip

If $G$ is an interpolation group with an order unit $u$, 
then $G$ is a unigroup with unit $u$, and the unit interval $E$ 
in $G$ has the Riesz-decomposition property \cite{Rav96}. In 
particular, if $G$ is a lattice-ordered abelian group with 
order unit $u$, then $G$ is a unigroup with unit $u$, and the 
unit interval $E$ in $G$ is an MV-algebra. If $V$ is a 
partially ordered vector space over any subfield of the 
real numbers and $u$ is an order unit in $V$, then, regarded 
as a partially-ordered additive abelian group, $V$ is a unigroup 
with unit $u$ \cite[Corollary 4.6]{BFInt}. In particular, with 
the identity operator $\mathbf{1}$ as order unit, the partially 
ordered additive group $\mathbb{G}({\cal H})$ of bounded 
self-adjoint operators on a Hilbert space ${\cal H}$ is a 
unigroup.

If an effect algebra $E$ can be realized as (or is isomorphic 
to) the unit interval in a unigroup $G$, then by definition  
$E$ is an IEA. Conversely, by \cite[Corollary 4.2]{BFInt} 
or by \cite[Theorem 5.4 and ff.]{FGB96}, {\em every IEA $E$ 
can be realized as the unit interval in a unigroup $G$}. 
Furthermore, {\em $G$ is uniquely determined by $E$ up to a 
unital isomorphism}, hence (by a slight abuse of language) 
we shall refer to $G$ as {\em the unigroup for} $E$. Thus, 
with $E$ as the unit interval in $G$, and $G$ as the unigroup 
for $E$, {\em we have a correspondence} $E\leftrightarrow G$ 
(up to isomorphism) {\em between IEA's $E$ and unigroups $G$}.  
More formally, {\em there is a categorical equivalence between 
the category of interval effect algebras and the category of 
unigroups} \cite[Theorem 3]{PEAwRD}.

If $G\not=\{0\}$ is a partially-ordered abelian group with order 
unit $u$, then a {\em state} on $G$ is defined to be a 
homomorphism $\omega\colon G\to\reals$ from $G$ to the additive 
group of real numbers such that $\omega(G\sp{+})\subseteq\reals
\sp{+}$ and $\omega(u)=1$ \cite[Chapter 4]{Good}. Denote by 
$\Omega(G)$ the set of all states on $G$. Then $\Omega(G)$ is a 
subset of the locally convex linear topological space $\reals
\sp{G}$ of all functions from $G$ to $\reals$ with pointwise 
operations and the topology of pointwise convergence. As such, 
$\Omega(G)$ is nonempty \cite[Corollary 4.4]{Good} and it is a 
compact convex subset of $\reals\sp{G}$ \cite[Proposition 6.2]
{Good}. Therefore, by the Krein-Milman theorem, $\Omega(G)$ is 
the closed convex hull of its own set $\partial\sb{e}\Omega(G)$ 
of extreme points.  By \cite[Theorem 4.14]{Good}, $G$ is 
archimedean iff $\Omega(G)$ determines $G\sp{+}$ in the sense 
that $G\sp{+}=\{g\in G\mid 0\leq\omega(g)\text{ for all }\omega
\in\Omega(G)\}$.

Suppose that $E$ is the unit interval in a unigroup $G\not=
\{0\}$.  If $\omega\in\Omega(G)$, then the restriction $\pi :=
\omega|\sb{E}$ of $\omega$ to $E$ is a probability measure on 
$E$. Conversely, if $\pi\in\Pi(E)$, then $\pi\colon E\to\reals$ 
is an $\reals$-valued measure, hence it admits a unique 
extension to a homomorphism $\omega\colon G\to\reals$ into 
the additive group of real numbers. Moreover, as $\pi(E)
\subseteq\reals\sp{+}$ and $E$ generates $G\sp{+}$, it follows 
that $\omega\in\Omega(G)$. Thus, we have an affine isomorphism 
$\pi\leftrightarrow\omega$ with $\pi=\omega|\sb{E}$ between the 
space $\Pi(E)$ of probability measures on $E$ and the state 
space $\Omega(G)$ of $G$. As a consequence, if $G$ is 
archimedean, then $\Pi(E)$ is an order-determining set of 
probability measures on $E$.

Let $E$ be an IEA and let $G$ be the unigroup with unit $u$ 
for $E$. By \cite{Rav96}, $E$ has the Riesz-decomposition 
property iff $G$, the {\em Ravindran group} of $E$, is an 
interpolation group. By \cite{Mund}, $E$ is an MV-effect 
algebra iff $G$, the {\em Mundici group} of $E$, is lattice 
ordered. If $E$ is totally ordered, then so is its Mundici 
group $G$ \cite[Corollary 6.5]{BFInt}. By definition, $G$ is 
a {\em Boolean unigroup} iff $E$ is a Boolean effect algebra. 
By \cite[Theorem 4.26]{Ritt}, $G$ is a Boolean unigroup iff 
it is an interpolation group and $u$ is the smallest order unit 
in $G$. In Theorem 8.7 below, we characterize $G$ for the case 
in which $E$ is an HMV-effect algebra. In this connection, the 
following theorem of N. Ritter \cite[Theorem 4.21]{Ritt} is of 
interest (see Example 8.8 below).

\medskip

\noindent{\bf 7.4 Theorem} {\em Let $A$ be a lattice-ordered 
group with order unit $u$ and suppose that the $u$-interval  
in $A$ is an HMV-effect algebra. Then, if $0\not=v\in A\sp{+}$, 
the $v$-interval in $A$ is also an HMV-effect algebra.}
  
\section{Compressions} 

Let ${\cal B}$ be a field of subsets of the nonempty set $X$ and 
let ${\cal F}({\cal B},\integers)$ be the commutative  
lattice-ordered ring with unit defined in Example 1.1. Then,  
regarded as an additive lattice-ordered group with order unit 
$1$, $G :={\cal F}({\cal B},\integers)$ is a Boolean unigroup and 
$E :=G\sp{+}[0,1]$ is a Boolean effect algebra. The Boolean sum 
$x\oplus y$ on $E$ extends to the addition operation $g+h$ on 
$G$.  There are two natural options for extending the Boolean 
product $xy=x\wedge y$ on $E$ to $G$, namely (1) to the product 
operation $gh$, or (2) to the infimum operation $g\wedge h$, 
for $g,h\in G$.  Option (1), which Boole might have favored, 
can be generalized, but only to unigroups admitting a reasonable 
notion of a product (see \cite{FCOM}). Option (2), which Jevons 
might have preferred, can be generalized, but only to MV-algebras 
and their lattice-ordered unigroups.     

With the notation of the last paragraph, the key to a more 
general solution of the extension problem for Boolean-type 
products is as follows: Rather than looking at the product 
$xy$ as a binary operation, we fix $x\in E$ and consider the 
{\em unary operation} $y\mapsto xy$ for all $y\in E$. This unary 
operation has a natural extension to a unary operation $J\sb{x}
\colon G\to G$ defined by $J\sb{x}(g) :=xg$ for all $g\in G$. 
(One could imagine that Boole would have been comfortable with 
$J\sb{x}$ because of his work with differential operators.)
Evidently, $J\sb{x}$ is a ``retraction" on $G$ with ``focus" 
$x$ as per the following definition.

\medskip

\noindent{\bf 8.1 Definition} Let $G$ be a unital group with 
unit $u$ and unit interval $E$. A mapping $J\colon G\to G$ is 
called a {\em retraction} with {\em focus} $p$ iff $J$ is an 
order-preserving endomorphism on $G$ such that $p :=J(u)
\leq u$ and, for all $e\in E$, $e\leq p\Rightarrow J(e)=e$. 
A {\em compression} on $G$ is a retraction $J$ on $G$ such 
that, if $p$ is the focus of $J$, $e\in E$, and $J(e)=0$, 
then $e\leq p\sp{\perp}$.  A retraction $J\,'$ on $G$ is 
a {\em quasicomplement} of the retraction $J$ on $G$ iff, 
for all $g\in G\sp{+}$, $J(g)=0\Leftrightarrow J\,'(g)=g$ 
and $J\,'(g)=0\Leftrightarrow J(g)=g$.

\medskip

Every retraction $J$ on a unital group $G$ is idempotent, 
i.e., $J=J\circ J$ \cite[Lemma 2.2]{FCOM}. If a retraction 
$J$ has a quasicomplement $J\,'$, then both $J$ and $J\,'$ 
are compressions \cite[Lemma 3.2]{FCOM}.

\medskip

\noindent{\bf 8.2 Definition} A {\em compressible group} is a 
unital group $G$ such that every retraction on $G$ is determined 
by its focus and every retraction on $G$ has a quasicomplementary 
retraction on $G$.  Let $G$ be a compressible group with unit  
interval $E$. An element $p\in E$ is called a {\em projection} 
iff it is the focus of a retraction (hence a compression) on $G$. 
The set of all projections in $E$ is denoted by $P(G)$, and if 
$p\in P(G)$, the unique compression on $G$ with focus $p$ is 
denoted by $J\sb{p}$. 

\medskip

See \cite{FCG, FCOM, FCGGC} for the basic theory of compressible 
groups. Let $G$ be a compressible group with unit interval $E$.  
Then $P(G)$ is a subeffect algebra of $E$ and, as such, $P(G)$ is 
an orthomodular poset. If $p,q\in P(G)$, then $pCq\Leftrightarrow 
J\sb{p}\circ J\sb{q}=J\sb{q}\circ J\sb{p}$. In fact, if $pCq$, 
then $p\wedge q$ exists in $E$, $p\wedge q\in P(G)$, and $J\sb{p}
\circ J\sb{q}=J\sb{q}\circ J\sb{p}=J\sb{p\wedge q}$. Furthermore, 
the quasicomplement of $J\sb{p}$ is $J\sb{p\sp{\perp}}$.

If $E$ is an IEA and $G$ is the unigroup for $E$, then $G$ is a 
compressible group iff $E$ is a compressible effect algebra in 
the sense of Gudder \cite{Gud}. The notion of a compressible 
group enables a {\em multiplicative} characterization of a 
Boolean effect algebra.  Indeed, {\em if $E$ is the unit interval 
in a unital group $G$, then $E$ is a Boolean effect algebra 
iff $G$ is a compressible group with $E=P(G)$} \cite[Theorem 5.5]
{FCG}. The following theorem applies to the Ravindran group $G$ 
of an effect algebra $E$ with the Riesz-decomposition property.

\medskip

\noindent{\bf 8.3 Theorem} {\em Let $G$ be an interpolation group 
with order unit $u$. Then} (i) {\em $G$ is a compressible unigroup 
with unit $u$ and unit interval $E :=G\sp{+}[0,u]$, and $P(G)=
\{p\in E\mid p\wedge p\sp{\perp}=0\}=C(E)$ is a Boolean effect 
algebra. Let $p\in P(E)=C(E)$, $H :=J\sb{p}(G)$, and $K :=J
\sb{p\sp{\perp}}(G)$. Then:} (ii) $e\in E\Rightarrow J\sb{p}(e)=
p\wedge e$. (iii) $g\in G\Rightarrow g=J\sb{p}(g)+J\sb{p\sp
{\perp}}(g)$. (iv) {\em $H$ and $K$ are subgroups of $G$, 
$H+K=G$, and $H\cap K=\{0\}$.} (v) {\em With the partial order 
induced from $G$, $H$ and $K$ are interpolation groups with order 
units $p$ and $p\sp{\perp}$, respectively.} (vi) {\em The mappings 
$J\sb{p}\colon G\to H$ and $J\sb{p\sp{\perp}}\colon G\to K$ provide 
a projective representation of $G$ as a direct product of $H$ and 
$K$ in the category of unigroups.} (vii) {\em If $G$ is lattice 
ordered, then $H$ and $K$ are sublattices of $G$, and for all 
$g,h\in G$, $J\sb{p}(g\vee h)=J\sb{p}(g)\vee J\sb{p}(h)$, and 
$J\sb{p}(g\wedge h)=J\sb{p}(g)\wedge J\sb{p}(h)$.}  

\medskip

\noindent{\em Proof} For (i)--(vi), see \cite[Theorem 3.5]{FCG} and 
\cite[pp.127--131]{Good}. Part (vii) follows from (vi).\hspace{\fill}
$\square$  

\medskip

The compression operators on a compressible group provide a 
generaliza\-tion---in the spirit of Boole---of the multiplicative 
structure of an MV-algebra. This generalization is applicable 
not only to MV-algebras, but also to a large class of quantum 
logics as per the following example.

\medskip

\noindent{\bf Example 8.4} Let ${\cal A}$ be a unital C$\sp
{\ast}$-algebra with unity $1$. Then the additive group 
$G({\cal A})$ of self-adjoint elements in ${\cal A}$, partially 
ordered with the positive cone $G({\cal A})\sp{+}=\{aa\sp{\ast}
\mid a\in{\cal A}\}$, is an archimedean compressible unigroup  
\cite[Corollary 4.6]{FCOM}. The projections in $P(G({\cal A}))$ 
are the idempotent elements $p=p\sp{2}=p\sp{\ast}\in G({\cal A})$, 
and for $a\in G({\cal A})$, $J\sb{p}(a)=pap$.\hspace{\fill}$
\square$

\medskip 

\noindent{\bf 8.5 Definition} Let $G$ be a compressible group, 
let $g\in G$, and $p\in P(G)$. Then (i) $C(p) :=
\{g\in G\mid g=J\sb{p}(g)+ J\sb{u-p}(g)\}$ and (ii) $CPC(g) 
:=\bigcap\{C(p)\mid p\in P(G)\text{ and }g\in C(p)\}$. Elements 
$g\in C(p)$ are said to be {\em compatible} with the projection 
$p$.

\medskip

It is not difficult to verify that the notion of compatibility 
in Definition 8.5 (i) is consistent with the notion of compatibility 
in Definition 3.10, i.e., if $E$ is the unit interval in the 
compressible group $G$, then for $e\in E$ and $p\in P(G)$, 
$e\in C(p)\Leftrightarrow eCp$. If $g\in G$, then $CPC(g)$ 
in Definition 8.5 (ii) is the set of all elements $h\in G$ that 
are compatible with every projection with which $g$ is compatible. 
By Theorem 8.3, if $G$ is an interpolation unigroup, then $G=
C(p)$ for every $p\in P(G)$, hence $G=CPC(g)$ for every $g\in G$.

In Example 8.4, if $g\in G({\cal A})$ and $p\in P(G({\cal A}))$, 
then $g\in C(p)$ iff $gp=pg$, i.e., iff $g$ commutes with $p$ in 
the C$\sp{\ast}$-algebra ${\cal A}$. Hence, if ${\cal A}$ is a von 
Neumann algebra and $g,h\in G({\cal A})$, then $h\in CPC(g)$ iff 
$h$ ``double commutes" with $g$, i.e., $h$ commutes with every 
element $a\in{\cal A}$ that commutes with $g$.

\medskip

\noindent{\bf 8.6 Definition} Let $G$ be a compressible group.
\begin{list}%
{(\roman{rom})}{\usecounter{rom}
\setlength{\rightmargin}{\leftmargin}}
\item $G$ has the {\em Rickart projection property} iff 
there exists a mapping $'\colon G\to P(G)$, called the {\em 
Rickart mapping}, such that, for all $g\in G$ and all $p\in 
P(G)$, $p\leq g\,'\Leftrightarrow g\in C(p)$ with $J\sb{p}(g)
=0$. 
\item $G$ has the {\em general comparability} (GC) {\em 
property} iff, for every $g\in G$, there exists $p\in P(G)$ 
such that $p\in CPC(g)$ and $J\sb{p\sp{\perp}}(g)\leq 0\leq 
J\sb{p}(g)$.
\item An {\em RGC-group} (also called a {\em Rickart comgroup} 
\cite{FSpec}) is a compressible group with both the Rickart 
projection and general comparability properties. 
\end{list}

With $1$ as the order unit, the group $G({\cal A})$ of 
self-adjoint elements in a von Neumann algebra ${\cal A}$ 
is an archimedean RGC-unigroup. In an RGC-group $G$, every 
element has a rational spectral resolution, which, if $G$ is  
archimedean, has the basic properties of the spectral 
resolution of a bounded self-adjoint operator on a Hilbert 
space \cite{FSpec}.

Let $G$ be an interpolation group with order unit $u$. Then the 
conditions $g\in C(p)$ in Definition 8.6 (i) and $p\in CPC(g)$ 
in Definition 8.6 (ii) hold automatically, and $G$ satisfies 
the general comparability property iff it satisfies the condition 
with the same name in \cite[Chapter 8]{Good}.  Thus, by 
\cite[Proposition 8.9]{Good}, if $G$ satisfies the GC-property, 
then it is lattice ordered, and by \cite[Proposition 9.9]{Good}, 
if $G$ is lattice-ordered and Dedekind $\sigma$-complete, then 
$G$ has the GC-property. As a consequence of the following 
theorem, {\em if $G$ is the Mundici group of an MV-algebra $E$, 
then $E$ is an HMV-algebra iff $G$ is an RGC-group.}

\medskip

\noindent{\bf 8.7 Theorem} {\em Let $G$ be a unital group 
with order unit $u$ and unit interval $E$.  Then the following 
conditions are mutually equivalent:} (i) {\em $G$ is a unigroup 
and $E$ is an HMV-algebra.} (ii) {\em $G$ is lattice ordered, 
has the Rickart projection property, and the Rickart mapping 
$g\mapsto g\,'$ satisfies $e\wedge f=0\Rightarrow e\leq f\,'$ 
for all $e,f\in E$.} (iii) {\em $G$ is an RGC-group and 
$P(G)\subseteq C(E)$.}

\medskip

\noindent{\em Proof} (i) $\Rightarrow$ (ii). Assume (i) and 
let $'\colon E\to E$ be the Heyting negation connective. As 
$G$ is a unigroup, it is the Mundici group of the MV-algebra 
$E$, hence $G$ is lattice ordered, so it is an interpolation 
group. By Theorem 8.3, $G$ is a compressible group and $P(G)=
C(E)$ is a Boolean algebra. Let $g\in G\sp{+}$. Then there are 
elements $e\sb{1},e\sb{2},...,e\sb{n}\in E$ such that $g=\sum
\sb{i=1}\sp{n}e\sb{i}$. Let $q :=(e\sb{1})'\wedge(e\sb{2})'
\wedge\cdots\wedge(e\sb{n})'\in P(G)=C(E)$ and let $p\in P(G)$. 
Then $J\sb{p}(g)=0\Leftrightarrow\sum\sb{i=1}\sp{n}J\sb{p}
(e\sb{i})=0$, and since $0\leq J\sb{p}(e\sb{i})$ for $i=
1,2,...,n$, it follows that $J\sb{p}(g)=0\Leftrightarrow 
J\sb{p}(e\sb{i})=0\text{ for }i=1,2,...,n.$ But, by Theorem 
8.3, $J\sb{p}(e\sb{i})=0\Leftrightarrow p\wedge e\sb{i}=0
\Leftrightarrow p\leq (e\sb{i})'$ for $i=1,2,...,n$, and it 
follows that $J\sb{p}(g)=0\Leftrightarrow p\leq q$. Thus,  
$q$ depends only on $g$ and is independent of the choice of 
$e\sb{1},e\sb{2},...,e\sb{n}\in E$, hence we can and do extend 
the Heyting negation $'$ from $E$ to $G\sp{+}$ by defining $g\,' 
:=q\in P(G)$.

If $g\in G$, then $g\in G\sp{+}\Leftrightarrow g=g\sp{+}
\text{ and }g\sp{-}=0$. Therefore, we can and do extend $'$ 
from $G\sp{+}$ to $G$ by defining $g\,' :=(g\sp{+})'\wedge
(g\sp{-})'$ for all $g\in G$. Let $g\in G$ and let $p\in 
P(G)=C(E)$. By Theorem 8.3 (vii), $J\sb{p}(g\sp{+})=J\sb{p}
(g\vee 0)=J\sb{p}(g)\vee J\sb{p}(0)=J\sb{p}(g)\vee 0=
(J\sb{p}(g))\sp{+}$.  Thus, $J\sb{p}(g\sp{-})=J\sb{p}
((-g)\sp{+})=(J\sb{p}(-g))\sp{+}=(-J\sb{p}(g))\sp{+}=
(J\sb{p}(g))\sp{-}$, and it follows that $J\sb{p}(g)=0
\Rightarrow J\sb{p}(g\sp{+})=J\sb{p}(g\sp{-})=0\Rightarrow 
p\leq(g\sp{+})'\wedge(g\sp{-})'=g\,'.$  Conversely, $p\leq g\,'
\Rightarrow J\sb{p}(g\sp{+})=J\sb{p}(g\sp{-})=0\Rightarrow 
J\sb{p}(g)=J\sb{p}(g\sp{+}-g\sp{-})=J\sb{p}(g\sp{+})-
J\sb{p}(g\sp{-})=0$.  Therefore, $G$ has the Rickart 
projection property, and, since the Rickart mapping $'$ 
is an extension of the Heyting negation connective on 
$E$, it satisfies $e\wedge f=0\Rightarrow e\leq f\,'$ 
for all $e,f\in E$.

(ii) $\Rightarrow$ (iii). Assume (ii) and let $'\colon G
\to P(G)$ be the Rickart mapping on $G$. As $G$ is lattice 
ordered, it has the interpolation property and is a 
compressible group with $P(G)=C(E)$. By \cite[Lemma 6.2 
(iii)]{FCGGC}, $p\,'=p\sp{\perp}$ for all $p\in P(G)$. 
 
Suppose $a,b\in G\sp{+}$ and $a\wedge b=0$. Then there are 
elements $e\sb{i}\in E$ for $i=1,2,...,n$ and $f\sb{j}\in E$ 
for $j=1,2,...,m$ such that $a=\sum\sb{i=1}\sp{n}e\sb{i}$ 
and $b=\sum\sb{j=1}\sp{m}f\sb{j}$, and since $e\sb{i}\wedge 
f\sb{j}\leq a\wedge b=0$, it follows from (ii) that $e\sb{i}
\leq(f\sb{j})\,'$ for $i=1,2,...,n$ and $j=1,2,...,m$. By 
\cite[Lemma 6.2 (viii)]{FCGGC}, it follows that $(e\sb{i})
\,''\leq(f\sb{j})\,'$ for $i=1,2,...,n$ and $j=1,2,...,m$. 
By \cite[Theorem 6.4 (iv)]{FCGGC}, $a\,''=\bigvee\sb{i=1}
\sp{n}(e\sb{i})\,''$ and $b\,'=\bigwedge\sb{j=1}\sp{m}
(f\sb{j})\,'$, whence $a\,''\leq b\,'$.  Therefore, if 
$a,b\in G\sp{+}$, then $a\wedge b=0\Rightarrow a\,''\leq 
b\,'$.

Let $g\in G$ and define $q :=(g\sp{+})\,'$ and $p :=
(g\sp{-})\,'$. Then, as $g\sp{+},g\sp{-}\in G\sp{+}$ and 
$g\sp{+}\wedge g\sp{-}=0$, it follows that $p\sp{\perp}=
p\,'=(g\sp{-})\,''\leq (g\sp{+})\,'=q$.  Since $0=J\sb{q}
(g\sp{+})=J\sb{q}(g\vee 0)=J\sb{q}(g)\vee 0$, it follows 
that $J\sb{q}(g)\leq 0$.  Likewise, $J\sb{p}(-g)\leq 0$, 
i.e., $0\leq J\sb{p}(g)$. As $p\sp{\perp}\leq q$, we 
have $J\sb{p\sp{\perp}}(g)=J\sb{p\sp{\perp}\wedge q}(g)
=J\sb{p\sp{\perp}}(J\sb{q}(g))\leq 0$.  Therefore, $J\sb
{p\sp{\perp}}(g)\leq 0\leq J\sb{p}(g)$, so $G$ has the 
GC-property, and hence it is an RGC-group.  

(iii) $\Rightarrow$ (i). Assume (iii).  Let $p\in P(G)
\subseteq C(E)$. Then, if $e\in E$, we have $pCe$, hence 
$e\in C(p)$. If $g\in G$, then $g$ is finite linear 
combination with integer coefficients of elements of 
$E$, hence $g\in C(p)$. As $G$ has the GC-property and 
$g\in C(p)$ for all $g\in G$, $p\in P(G)$, it follows from 
\cite[Theorem 4.9]{FCG} that $G$ is lattice ordered. 
Therefore $G$ is a unigroup, $E$ is an MV-algebra, $G$ is 
the Mundici group of $E$, and $P(G)=C(E)$. 

Suppose $e,f\in E$ with $e\wedge f=0$. By the GC-property, 
there is a projection $p\in P(G)=C(E)$ such that $J\sb
{p\sp{\perp}}(e-f)\leq 0\leq J\sb{p}(e-f)$, i.e., $p\sp
{\perp}\wedge e=J\sb{p\sp{\perp}}(e)\leq J\sb{p\sp{\perp}}
(f)=p\sp{\perp}\wedge f$ and $p\wedge f=J\sb{p}(f)\leq 
J\sb{p}(e)=p\wedge e$. Thus, $J\sb{p\sp{\perp}}(e)=p\sp
{\perp}\wedge e=p\sp{\perp}\wedge e\wedge f=0$ and $J\sb{p}
(f)=p\wedge f=p\wedge f\wedge e=0$, and it follows that 
$p\,'=p\sp{\perp}\leq e\,'$ and $p\leq f\,'$. Therefore, 
$e\leq e\,''\leq p\,''=p\leq f\,'$, and we conclude that,  
for $e,f\in E$, $e\wedge f=0\Rightarrow e\leq f\,'$. Conversely, 
suppose $e,f\in E$ and $e\leq f\,'$. Then, by \cite[Lemma 6.2 
(vii)]{FCGGC}, $e\,''\leq f\,'$, so $e\,''\wedge f=J\sb{e\,''}
(f)=0$, whence by \cite[Lemma 6.2 (vii)]{FCGGC}, $0\leq e\wedge 
f\leq e\,''\wedge f=0$. Consequently, the restriction to $E$ of 
$'$ satisfies the condition $e\wedge f=0\Leftrightarrow e\leq 
f\,'$ for all $e,f\in E$.  Therefore, by \cite[Theorem 6.13]{FGB98}, 
$E$ is a Heyting algebra under the Heyting conditional $(e
\supset f) :=(e-e\wedge f)\,'\vee f=((e-f)\sp{+})\,'\vee f$ for 
$e,f\in E$.\hspace{\fill}$\square$

\medskip

The following example, which generalizes Example 1.1, provides 
a large class of archimedean lattice-ordered RGC-unigroups with 
HMV-algebras as their unit intervals.

\medskip

\noindent{\bf 8.8 Example} Let ${\cal{B}}$ be a field of 
subsets of a set $X$, let $A$ be a lattice-ordered unigroup 
such that the unit interval in $A$ is an HMV-effect algebra 
(e.g., $A=\integers$ with unit $1$). Define ${\cal{F}}
({\cal{B}},A)$ to be the partially ordered abelian group 
under pointwise addition and pointwise partial order of all 
functions $f\colon X\to A$ such that (i) $f\sp{-1}(a)\in
{\cal{B}}$ for all $a\in A$ and (ii) $\{f(x)\mid x\in X\}$ is 
a finite subset of $A$. Then an element $u\in{\cal{F}}
({\cal{B}},A)$ is an order unit iff $u(x)$ is an order unit 
in $A$ for all $x\in X$. If $u$ is an order unit in ${\cal{F}}
({\cal{B}},A)$, then ${\cal{F}}({\cal{B}},A)$ is a 
lattice-ordered RGC-unigroup with unit $u$, hence the 
u-interval $E$ in ${\cal{F}}({\cal{B}},A)$ is an HMV-algebra.
\hspace{\fill}$\square$

\medskip

\noindent{\bf Acknowledgment} The author is grateful for 
stimulating conversations with Professor Richard J. Greechie 
regarding the subject matter of this article.  However, any 
errors or misapprehensions are strictly the fault of the 
author.


\begin{thebibliography}{99}


\bibitem{BFPhi} Bennett, M.K. and Foulis, D.J., Phi-symmetric effect 
algebras, {\em Foundations of Physics} {\bf 25}, No. 12 (1995) 
1699--1722.

\bibitem{BFInt} Bennett, M.K. and Foulis, D.J.,  Interval and scale     
  effect algebras, {\em Advances in Applied Mathematics} {\bf 19} 
  (1997) 200--215.

\bibitem{Beran} Beran, L., {\em Orthomodular Lattices, An Algebraic 
Approach}, Mathematics and its Applications, Vol. 18, D. Reidel Publishing 
Company, Dordrecht, 1985.

\bibitem{Boo} Boole, G., {\em The Laws of Thought}, Dover Publications, 
Inc. N.Y., first printing 1854.

\bibitem{BLM} Busch, P., Lahti, P.J., and Mittelstaedt, P., {\em The 
Quantum Theory of Measurement}, Lecture Notes on Physics, Vol. 2, 
Springer, Berlin. (1991).

\bibitem{CCC} Chang, C.C., Algebraic analysis of many-valued 
logics, {\em Transactions of the American Mathematical Society}
{\bf 88} (1957) 467--490.

\bibitem{CK} Chovanec, F. and K\^{o}pka, F., Boolean D-posets, {\em 
Tatra Mountains Mathematical Publications} {\bf 10} (1997) 183--197.

\bibitem{DG89} Dalla Chiara, M.L. and Giuntini, R., Paraconsistent 
quantum logics, {\em Foundations of Physics} {\bf 19}, No. 7 
(1989) 891--904. 

\bibitem{DGG} Dalla Chiara, M.L., Giuntini, R., and Greechie, R.J., 
{\em Reasoning in Quantum Theory}, {\em Trends in Logic}, Vol 22, 
Kluwer, Dordrecht/Boston/London, 2004. 

\bibitem{Fou62} Foulis, D.J.,  A note on orthomodular lattices, 
{\em Portugaliae Mathematica} {\bf 21} (1962) 65--72.

\bibitem{Fou89} Foulis, D.J., Coupled physical systems, {\em 
Foundations of Physics}, {\bf 7} (1989) 905--922.

\bibitem{FouMV} Foulis, D.J., MV and Heyting effect algebras, 
{\em Foundations of Physics} {\bf 30}, No. 10 (2000) 1687--1706.

\bibitem{FCG} Foulis, D.J., Compressible groups, 
{\em Mathematica Slovaca} {\bf 53}, No. 5 (2003)
433--455.
 
\bibitem{FCOM} Foulis, D.J., Compressions on partially 
ordered abelian groups, to appear in {\em Proceedings 
of the American Mathematical Society}.

\bibitem{FCGGC} Foulis, D.J., Compressible groups with 
general comparability, to appear in {\em Mathematica 
Slovaca}.

\bibitem{FSpec} Foulis, D.J., Spectral resolution in a Rickart 
comgroup, to appear in {\em Reports on Mathematical Physics}.

\bibitem{FB} Foulis, D.J. and Bennett, M.K., Effect algebras and 
unsharp quantum logics,  {\em Foundations of Physics} {\bf 24}, 
No. 10 (1994) 1325--1346.

\bibitem{FGB96} Foulis, D.J., Greechie, R.J., and Bennett, 
M.K., Test groups and effect algebras, {\em International 
Journal of Theoretical Physics} {\bf 35}, No. 6 (1996) 
1117--1140.

\bibitem{FGB98} Foulis, D.J., Greechie, R.J., and Bennett, M.K., 
The transition to unigroups, {\em International Journal of 
Theoretical Physics} {\bf 37}, No. 1 (1998) 45--63.

\bibitem{Good}  Goodearl, K.R., {\em Partially Ordered 
Abelian Groups with Interpolation}, A.M.S. Mathematical 
Surveys and Monographs, {\bf No. 20}, American 
Mathematical Society, Providence, RI, 1986.

\bibitem{GFP} Greechie, R.J., Foulis, D.J., and Pulmannov\'{a} S.,  
The center of an effect algebra, {\em Order} {\bf 12} (1995) 
91--106.

\bibitem{Gud65}  Gudder, S.P., Spectral methods for a generalized 
probability theory, {\em Transactions of the American Mathematical 
Society} {\bf 119} (1965) 428--442.

\bibitem{Gud68}  Gudder, S.P., {\em Quantum Probability}, Academic 
Press, San Diego, 1988.    

\bibitem{Gud} Gudder, S.P., Compressible effect algebras, to 
appear.

\bibitem{IH76} Hailperin, I., {\em Boole's Logic and Probability, 
Studies in Logic and the Foundations of Mathematics} {\bf 85} 
North-Holland Publishing Co., Amsterdam/New York/Oxford, 1976

\bibitem{IH81} Hailperin, I., Boole's algebra isn't Boolean 
algebra, {\em Mathematics Magazine} {\bf 54}, No. 4 (1981) 
172--184.

\bibitem{Kalm} Kalmbach, G., {\em Orthomodular Lattices}, 
Academic Press, Inc., London/New York, 1983. 

\bibitem{Mund} Mundici, D., Interpretation of AF C$\sp
{\ast}$-algebras in \L ukasiewicz sentential calculus, 
{\em Journal of Functional Analysis} {\bf 65} (1986) 15--63.

\bibitem{JvN} Neumann, J. von {\em Grundlagen der Quantenmechanik}, 
Springer Verlag, Berlin/Heidelberg/New York, 1932.  English 
translation: {\em Mathematical Foundations of Quantum mechanics}
Princeton University Press, Princeton, NJ, 1955.

\bibitem{PP}Pt\'{a}k, P. and Pulmannov\'{a}, S., {\em 
Orthomodular Structures as Quantum Logics}, Kluwer, 
Dordrecht/Boston/London, 1991.

\bibitem{PEAwRD} Pulmannov\'{a}, S., Effect algebras with 
the Riesz decomposition property and AF C$\sp{\ast}$-algebras, 
{\em Foundations of Physics} {\bf 29}, No. 9 (1999)1389--1401.

\bibitem{JP} Pykacz, J., Conjunctions, disjunctions, and 
Bell-type inequalities in orthoalgebras, {\em International 
Journal of Theoretical Physics} {\bf 35}, No. 11 (1996) 
2353--2363.

\bibitem{Rav96} Ravindran, K., {\em On a Structure Theory of Effect  
Algebras}, Ph.D. Thesis, Kansas State University, 1996, UMI 
Dissertation Services, No. 9629062, Ann Arbor, MI.

\bibitem{Ritt} Ritter, N., {\em Order Unit Intervals in Unigroups},  
Ph.D. Thesis, University of Massachusetts, Department of Mathematics 
and Statistics, 2000. 

\end{thebibliography}
\end{document}